\documentclass{ifacconf}

\usepackage{graphicx}      
\usepackage{natbib}        
\usepackage{amsmath} 
\usepackage{amssymb}  
\usepackage{graphicx}
\usepackage{amsmath}
\usepackage{epstopdf}
\usepackage{mathtools}
\usepackage{dsfont}
\usepackage{balance}

\DeclareMathOperator{\Succ}{Succ}
\DeclareMathOperator{\Pre}{Pre}
\DeclareMathOperator{\Int}{Int}
\newcommand{\beq}{\begin{equation}}
	\newcommand{\eeq}{\end{equation}}
\newcommand{\rr}{{\mathbb R}}
\newcommand{\zz}{{\mathbb Z}}
\newcommand{\cvd}{\hfill$\blacksquare$\vspace{0em}\par\noindent}

\newcounter{algorithmctr}[section]
\renewcommand{\thealgorithmctr}{\thesection.\arabic{algorithmctr}}
{\refstepcounter{algorithmctr}\begin{list}{}{%
			\setlength{\rightmargin}{0\linewidth}%
			\setlength{\leftmargin}{.05\linewidth}
			\setlength{\itemsep}{1pt}
			\setlength{\parskip}{0pt}
			\setlength{\parsep}{0pt}}%
		\rmfamily\small
		\item[]{\setlength{\parskip}{0ex}\hrulefill\par%
			\nopagebreak{\bfseries\textsf{Algorithm \thealgorithmctr~}}}}%
	{{\setlength{\parskip}{-1ex}\nopagebreak\par\hrulefill} \end{list}}
\newtheorem{assumption}{Assumption}
\newtheorem{theorem}{Theorem}
\newtheorem{remark}{Remark}

\usepackage{amsmath} 
\usepackage{amssymb}  
\usepackage{graphicx}
\usepackage{amsmath}
\usepackage{epstopdf}
\usepackage{mathtools}
\usepackage{dsfont}
\begin{document}
\begin{frontmatter}

\title{Learning Model Predictive Control for Iterative Tasks: A Computationally Efficient Approach for Linear System} 

\author[First]{Ugo Rosolia} 
\author[First]{Francesco Borrelli} 

\address[First]{University of California at Berkeley ,
	Berkeley, CA 94701, USA (e-mail: \{ugo.rosolia, fborrelli\}@berkeley.edu).}

\begin{abstract}                
A Learning Model Predictive Controller (LMPC) for linear system is presented. 
The proposed controller builds on previous work on nonlinear LMPC and decreases its computational burden for linear system. The control scheme is reference-free and is able to improve its performance by learning from previous iterations. A convex safe set and a terminal cost function are used in order to guarantee recursive feasibility and non-increasing performance at each iteration.
The paper presents the control design approach, and shows how to recursively construct the convex terminal set and the terminal cost from state and input trajectories of previous iterations.
Simulation results show the effectiveness of the proposed control logic.
\end{abstract}


\end{frontmatter}
\section{INTRODUCTION}
Iterative Learning Control (ILC) studies control design for autonomous systems performing repetitive tasks \cite{c3,c7,c8}. One task execution is often referred to as ``iteration" or ``trial". In ILC, at each iteration, the system starts from the same initial condition and the controller objective is to track a given reference, rejecting periodic disturbances \cite{c3,c7}. The tracking error from the previous iterations is used to improve the tracking performance of the closed loop system. Different strategies have been proposed to guarantee zero tracking error of the closed loop system \cite{c3, c7, c8}.

Several control frameworks which combine ILC and MPC strategies have been proposed in literature, \cite{subbaraman2016extremum, lee2000convergence, lee2000model}.
In the classical ILC approach the goal of the controller is to track a reference trajectory, however, in some application such has autonomous racing \cite{c17,c18} or for some manipulation tasks \cite{2016learning}, it may be challenging to generate a priori a reference trajectory that maximize the system performance. For this reason, a very recent work \cite{2016learning} proposed a reference-free ILC  scheme. The authors used a MPC controller with a terminal cost that allows to consider the long term planning. This terminal cost is computed using a neural network trained on data generated by offline simulations. The authors were able to improve the system performance over iterations. However, no guaranties about stability, recursive feasibility and performance improvement are provided.

Our objective is to design a reference-free iterative control strategy for linear system able to learn from previous iterations. At each iteration,  the initial condition, the constraints and the objective
function do not change. The $j$-th iteration  cost is defined as the objective function evaluated for the realized closed loop system trajectory. The iteration cost shall not increase over the iterations and  state and input constraints shall be satisfied. Model Predictive Control is an appealing technique to tackle this problem for its ability to handle state and inputs constraints while minimizing a finite-time predicted cost \cite{c11}. However, the receding horizon nature can lead to infeasibility and it does not guaranty improved  performance at each iteration \cite{c12}.

The contribution of this paper is the following. We present an extension to  the learning MPC for iterative control task in \cite{LMPC}. In particular,  we introduce a new formulation for linear system that drastically reduces the computation burden of the controller without compromising the guaranties of the learning MPC. 
We show how to design a convex safe set and a terminal cost function
in order to guarantee: \emph{(i):} [\textit{asymptotic stability}], the closed loop system converges asymptotically to the equilibrium point, \emph{(ii):} [\textit{persistent feasibility}], state and input constraints are satisfied if they were satisfied at iterations $j-1$ \emph{(iii):} [\textit{performance improvement}], the $j$-th iteration cost does not increase compared with the $j$-$1$-th iteration cost, \emph{(iv):} [\textit{global optimality}], if the steady state system converges to a closed-loop trajectory as the number of iterations $j$ goes to infinity, then that closed-loop trajectory is globally optimal. We emphasize that \emph{(i)}-\emph{(ii)} are standard MPC design requirement and \emph{(iii)}-\emph{(iv)} are the core contribution of this work. 

This paper is organized as follows: in Section II we introduce the notation used throughout the paper. Then, we define the convex safe set and the terminal cost function used in the design of the learning MPC. Section III describes the control design. We show the recursive feasibility and stability of the control logic and, afterwards, we prove the convergence properties.
Finally, in Section IV we test the proposed control logic on an infinite horizon linear quadratic regulator and we compare the computational efficiency with the learning MPC from \cite{LMPC}.

\section{PROBLEM FORMULATION}

Consider the discrete time system
\begin{equation}\label{eq:system}
	x_{t+1}=A x_t + B u_t,
\end{equation}
where $x\in \rr^n$ and $u\in\rr^{m}$ are the system state and input,
respectively, subject to  the constraints
\begin{equation}
	x_t \in \mathcal{X},\ u_t\in \mathcal{U},\ \forall t\in\zz_{0+}.
	\label{eq:inv_constraints}
\end{equation}
where $\mathcal{X}$ and $\mathcal{U}$ are convex sets.

At the $j$-th iteration the vectors
\begin{subequations}\label{eq:sequence}
	\begin{align}
		{\bf{u}}^j ~ = ~ [u_0^j,~u_1^j,~...,~u_t^j,~...], \label{eq:sequenceU} \\
		{\bf{x}}^j ~ = ~ [x_0^j,~x_1^j,~...,~x_t^j,~...], \label{eq:sequenceX}
	\end{align}
\end{subequations}
collect the inputs applied to system~(\ref{eq:system}) and the corresponding
state evolution. In (\ref{eq:sequence}), $x_t^j$ and $u_t^j$ denote the system state and the control input at time $t$ of the $j$-th iteration, respectively. We assume that at each $j$-th iteration the closed loop trajectories start from the same initial state,
\begin{equation}\label{eq:Initial_state}
	\begin{aligned}
		x_0^j ~ = x_S, ~\forall j \geq 0.\\
	\end{aligned}
\end{equation}

The goal is to design a controller which solves the following infinite horizon optimal control problem at each iteration:
\begin{subequations}\label{eq:ConstraintsInf}
	\begin{align}
		J_{0\rightarrow \infty}^*(x_S)&=\min_{u_0,u_1,\ldots} \sum\limits_{k=0}^{\infty} h(x_k,u_k)\label{eq:Constraints1InfCost}\\
		\textrm{s.t. }
		&x_{k+1}=A x_k + B u_k,~\forall k\geq 0 \label{eq:Constraints1Inf}\\
		&x_0=x_S,\label{eq:Constraints2Inf}\\
		&x_k \in \mathcal{X},~u_k \in \mathcal{U},~\forall k\geq 0\label{eq:Constraints3Inf}
	\end{align}
\end{subequations}
where equations (\ref{eq:Constraints1Inf}) and (\ref{eq:Constraints2Inf}) represent the system dynamics and the initial condition, and (\ref{eq:Constraints3Inf}) are the state and input constraints. The stage cost, $h(\cdot,\cdot)$, in equation (\ref{eq:Constraints1InfCost}) is continuous, jointly convex and it satisfies
\begin{equation}\label{eq:RunningCost}
	\begin{aligned}
		h(x_F,0) = 0~\textrm{and}~ h(x_t^j,u_t^j) \succ 0 ~ \forall ~ x_t^j \in&~\rr^n \setminus \{x_F\},\\
		& u_t^j \in \rr^m\setminus \{0\},
	\end{aligned}
\end{equation}
where the final state $x_F$ is assumed to be a feasible equilibrium for the unforced system (\ref{eq:system})
\begin{equation}\label{eq:Def_x_f}
	\begin{aligned}
		x_F = A x_F.\\
	\end{aligned}
\end{equation}

Next we introduce the definition of the convex safe set and of the terminal cost. Both will be used later to guarantee stability and feasibility of the learning MPC for linear system.

\subsection{Convex Safe Set}
In the following we recall the definition of the sampled Safe Set from \cite{LMPC} which is necessary to construct the convex Safe Set used in the learning MPC for linear system.

The definition of the \emph{sampled Safe Set} exploits the iterative nature of the control task to define an invariant control set, using the realized system trajectories. At the $j$-th iteration the sampled safe set, $\mathcal{SS}^j$, is defined as
\begin{equation}\label{eq:SS}
	\begin{aligned}
		\mathcal{SS}^j = \textrm{}\left\{\bigcup_{i \in M^j} \bigcup_{t=0}^{\infty} x_t^i \right\}.
	\end{aligned}
\end{equation}
$\mathcal{SS}^j$ is the collection of all state trajectories at iteration $i$ for $i\in M^j$. $M^j$ in equation (\ref{eq:SS}) is the set of indexes $k$ associated with successful iterations $k$ for $k\leq j$, defined as:
\begin{equation}\label{eq:M}
	\begin{aligned}
		M^j = \textrm{} \Big\{ k \in [0,j] : \lim_{t \to \infty} x_t^k = x_F \Big\}.
	\end{aligned}
\end{equation}
Moreover, as $\mathcal{X}$ and $\mathcal{U}$ are convex, for each convex combination of the elements in  $\mathcal{SS}^j$ we can find a control sequence that steers the system (\ref{eq:system}) to $x_F$. Therefore, the \emph{convex Safe Set}, defined as
\begin{equation} \label{eq:CS}
	\begin{aligned}
		\mathcal{CS}^j = \text{Conv}(\mathcal{SS}^j) = \Big\{ \sum_{i=1}^{|\mathcal{SS}^j|} \alpha_i z_i : \alpha_i \geq 0, \sum_{i=1}^{|\mathcal{SS}^j|} \alpha_i = 1, \\z_i \in \mathcal{SS}^j \Big\},
	\end{aligned}
\end{equation}
is a control invariant set. Note that $|\mathcal{SS}^j|$ is the cardinality of $\mathcal{SS}^j$. For further details on control invariant set we refer to \cite{c16}.

From (\ref{eq:M}) we have that $M^i \subseteq M^j, \forall i \leq j$, which implies that
\begin{equation}\label{eq:CS_subset}
	\begin{aligned}
		\mathcal{CS}^i \subseteq \mathcal{CS}^j, \forall i \leq j.
	\end{aligned}
\end{equation}

\subsection{Terminal Cost}
At time $t$ of the $j$-th iteration the cost-to-go associated with the closed loop trajectory (\ref{eq:sequenceX}) and input sequence (\ref{eq:sequenceU}) is defined as

\begin{equation}\label{eq:Functional}
	\begin{aligned}
		J_{t\rightarrow \infty}^j(x_t^j) = ~ \sum\limits_{k=0}^{\infty} h(x_{t+k}^j,u_{t+k}^j),
	\end{aligned}
\end{equation}
where $h(\cdot,\cdot)$ is the stage cost of problem (\ref{eq:ConstraintsInf}).
We define the $j$-th iteration cost as the cost (\ref{eq:Functional}) of the $j$-th trajectory at time $t=0$,
\begin{equation}\label{eq:Performance}
	\begin{aligned}
		J_{0\rightarrow \infty}^j(x_0^j) = ~ \sum\limits_{k=0}^{\infty} h(x_k^j,u_k^j).
	\end{aligned}
\end{equation}
$J_{0\rightarrow \infty}^j(x_0^j)$ quantifies the controller performance at each $j$-th 
iteration.

\begin{remark}
	In equations (\ref{eq:Functional})-(\ref{eq:Performance}), $x_k^j$ and $u_k^j$ are the realized state and input at the $j$-th iteration, as defined in (\ref{eq:sequence}).
\end{remark}

Finally we define the, barycentric function (\cite{jones2010polytopic})
\begin{equation}
	P^j(x) = \begin{cases}
		p^{j,*}(x)  & \mbox{If } x \in \mathcal{CS}^j \\
		+\infty  & \mbox{If } x \notin \mathcal{CS}^j \\
	\end{cases}
\end{equation} 
where 
\begin{subequations} \label{p^*}
	\begin{align}
		p^{j,*}(x) &= \min_{\lambda_t^j \geq 0, \forall t \in [0, \infty)} \sum_{k=0}^{j} \sum_{t=0}^{\infty}\lambda_t^k J_{t\rightarrow \infty}^k(x_t^k)  \\
		s.t. & \notag \\
		&\sum_{k=0}^{j} \sum_{t=0}^{\infty}\lambda_t^k = 1 \\
		&\sum_{k=0}^{j} \sum_{t=0}^{\infty}\lambda_t^k x_t^k= x, \label{eq:p_term}
	\end{align}
\end{subequations}
where $x_t^k$ is the realized state at time $t$ of the $j$-th iteration, as defined in (\ref{eq:sequenceX}). \\

\begin{remark}
	The function $P^j(x)$ assigns to every point in $\mathcal{CS}^j$ the minimum cost-to-go along the trajectories in $\mathcal{CS}^j$, in particular we have that $\forall x \in \mathcal{CS}^j,$
	\begin{equation} \label{eq:Psol}
		\begin{aligned}
			P^j(x) =& \sum_{k =0}^{j} \sum_{t=0}^{\infty} \lambda^{*,k}_t J_{t\rightarrow \infty}^k(x_t^k)  = \\
			=& \sum_{k =0}^{j} \sum_{t=0}^{\infty} \lambda^{*,k}_t \sum\limits_{l=0}^{\infty} h(x_{t+l}^k,u_{t+l}^k) 
		\end{aligned}
	\end{equation}
	where $\lambda^{*,k}_t$ is the minimizer in (\ref{p^*}). \\
\end{remark}

\begin{remark}
	In practical applications each $j$-th iteration has a finite time duration $t_j$, and therefore $p^{j,*}(x)$ is reformulated as
	\begin{subequations} \label{eq:p^*Rel}
		\begin{align}
			p^{j,*}(x) &= \min_{\lambda_t \geq 0, \forall t \in [0, \infty)} \sum_{k=0}^{j} \sum_{t=0}^{t_j}\lambda_t^k J_{t\rightarrow \infty}^k(x_t^k) \\
			s.t. & \notag \\
			&\sum_{k=0}^{j} \sum_{t=0}^{t_j}\lambda_t^k = 1 \\
			&\sum_{k=0}^{j} \sum_{t=0}^{t_j}\lambda_t^k x_t^k= x.
		\end{align}
	\end{subequations}
\end{remark}

\section{LMPC FOR LINEAR SYSTEM}
In this section we present the
design of the proposed Learning Model Predictive Control (LMPC). We first assume that there exists a feasible input sequence that steers the system from the initial point $x_S$ to terminal point $x_F$ at the $0$-th iteration.
Then we prove that the proposed LMPC is guaranteed to be recursively feasible, i.e. feasible at all time instants of every successive iteration. Moreover, we show that the LMPC guaranties a non-increasing iterations cost between two successive executions of the task.

\subsection{LMPC Control Design}

The LMPC tries to compute a solution to the infinite time optimal control problem (\ref{eq:ConstraintsInf}) by solving at time $t$ of iteration $j$ the finite time constrained optimal control problem

\begin{subequations}\label{eq:Constraints}
	\begin{align}
		J_{t\rightarrow t+N}^{\scalebox{0.4}{LMPC},j}&(x_t^j)=\min_{u_{t|t},\ldots,u_{t+N-1|t}} \bigg[  \sum_{k=t}^{t+N-1}  h(x_{k|t},u_{k|t}) +\notag\\
		&~~~~~~~~~~~~~~~~~~~~~~~~~~~~~~~~~~~~~~~+ P^{j-1}(x_{t+N|t}) \bigg]\label{eq:LMPCobj}\\
		\textrm{s.t. }&\notag \\
		&x_{k+1|t}=A x_{k|t} + B u_{k|t},~\forall k \in [t, \cdots, t+N-1] \label{eq:Constraints1}\\
		&x_{t|t}=x_t^j,\label{eq:Constraints2}\\
		&x_{k|t} \in \mathcal{X}, ~ u_k \in \mathcal{U},~\forall k \in [t, \cdots, t+N-1] \label{eq:Constraints4}\\
		&x_{t+N|t} \in ~\mathcal{CS}^{j-1},\label{eq:Constraints5}
	\end{align}
\end{subequations}
where (\ref{eq:Constraints1}) and (\ref{eq:Constraints2}) represent the system dynamics and initial condition, respectively. The state and input constraints are given by (\ref{eq:Constraints4}). Finally (\ref{eq:Constraints5}) forces the terminal state into the set $\mathcal{CS}^{j-1}$ defined in equation (\ref{eq:CS}).\\
Let
\begin{equation}\label{eq:OptimalSolutionMPC}
	\begin{aligned}
		{\bf{u}}^{*,j}_{t:t+N|t}  &= [u_{t|t}^{*,j}, \cdots, u_{t+N-1|t}^{*,j}]\\
		{\bf{x}}^{*,j}_{t:t+N|t} &= [x_{t|t}^{*,j}, \cdots, x_{t+N|t}^{*,j}]
	\end{aligned}
\end{equation}
be the optimal solution of (\ref{eq:Constraints}) at time $t$ of the $j$-th iteration and $J_{t\rightarrow t+N}^{\scalebox{0.4}{LMPC},j}(x_t^j)$ the corresponding optimal cost. Then, at time $t$ of the iteration $j$, the first element of ${\bf{u}}^{*,j}_{t:t+N|t}$ is applied to the system (\ref{eq:system})
\begin{equation}\label{eq:MPC}
	\begin{aligned}
		u_t^j = u_{t|t}^{*,j}.
	\end{aligned}
\end{equation}
The finite time optimal control problem (\ref{eq:Constraints}) is repeated at time $t+1$, based on the new state $x_{t+1|t+1} = x_{t+1}^j$ (\ref{eq:Constraints2}), yielding a \textit{moving} or \textit{receding horizon} control strategy.

\begin{remark}
	Problem (\ref{eq:Constraints}) is a convex optimization problem as the terminal constraint (\ref{eq:LMPCobj})  enforces the terminal state in the convex set $\mathcal{CS}^{j-1}$ and the terminal cost $P^{j-1}(\cdot)$ in (\ref{eq:Constraints5}) is a convex function. This new formulation of the LMPC (\ref{eq:Constraints}), (\ref{eq:MPC}) as a convex problem is the main contribution of this work compared to \cite{LMPC}.
\end{remark}

\begin{assumption}
	\label{ass1}
	At iteration $j=1$ we assume that  $\mathcal{CS}^{j-1}= \mathcal{CS}^0$ is a non-empty set and that the trajectory ${\bf{x}}^0 \in \mathcal{CS}^0$ is feasible and convergent to $x_F$.
\end{assumption}

In the next section we prove that, under Assumption~\ref{ass1}, the LMPC (\ref{eq:Constraints}) and (\ref{eq:MPC}) in closed loop with system (\ref{eq:system}) guarantees recursively feasibility and stability, and non-increase of the iteration cost at each iteration.

\subsection{Recursive feasibility and stability}

In this Section, the properties of $\mathcal{CS}^j$ and $P^j(\cdot)$ are used to show recursive feasibility and asymptotic stability of the equilibrium point $x_F$.

\begin{theorem}
	Consider system (\ref{eq:system}) controlled by the \mbox{LMPC}
	controller (\ref{eq:Constraints}) and (\ref{eq:MPC}).
	Let $\mathcal{CS}^j$ be the convex safe set at iteration $j$ as defined in (\ref{eq:CS}). Let assumption 1 hold, then the LMPC (\ref{eq:Constraints}) and (\ref{eq:MPC}) is feasible $\forall~t \in\zz_{0+}$ and iteration $j\geq1$.
	Moreover, the equilibrium point $x_F$ is asymptotically stable for the closed loop system (\ref{eq:system}) and (\ref{eq:MPC}) at every iteration $j\geq1$.
\end{theorem}

\textbf{Proof:}
The proof follows from standard MPC arguments. By assumption $\mathcal{CS}^0$ is non empty. From (\ref{eq:CS_subset}) we have that $\mathcal{CS}^{0} \subseteq \mathcal{CS}^{j-1} ~ \forall j \geq 1$, and consequently $\mathcal{CS}^{j-1}$ is a non empty set. In particular, there exists a trajectory ${\bf{x}}^0 \in \mathcal{CS}^0 \subseteq \mathcal{CS}^{j-1}$. From (\ref{eq:Initial_state}) we know that $x_0^j = x_S~\forall j \ge 0$. At time $t=0$ of the $j$-th iteration the $N$ steps trajectory
\begin{equation}
	\begin{gathered}\label{eq:SolutionFeasible0}
		[x_0^0,~x_1^0,~...,x_{N}^0] \in \mathcal{CS}^{j-1},
	\end{gathered}
\end{equation}
and the related input sequence,
\begin{equation}
	\begin{gathered}\label{eq:SolutionFeasible}
		[u_0^0,~u_1^0,~...,u_{N-1}^0],
	\end{gathered}
\end{equation}
satisfy input and state constrains (\ref{eq:Constraints1})-(\ref{eq:Constraints2})-(\ref{eq:Constraints4}). Therefore (\ref{eq:SolutionFeasible0})-(\ref{eq:SolutionFeasible}) is a feasible solution to the LMPC (\ref{eq:Constraints}) and (\ref{eq:MPC}) at $t=0$ of the $j$-th iteration.\\
Assume that at time $t$ of the $j$-th iteration the LMPC (\ref{eq:Constraints}) and (\ref{eq:MPC}) is feasible and let ${\bf{x}}^{*,j}_{t:t+N|t}$ and ${\bf{u}}^{*,j}_{t:t+N|t}$
be the optimal trajectory and input sequence, as defined in (\ref{eq:OptimalSolutionMPC}). From (\ref{eq:Constraints2}) and (\ref{eq:MPC}) the realized state and input at time $t$ of the $j$-th iteration are given by
\begin{equation}\label{eq:BC}
	\begin{aligned}
		x_t^j = x_{t|t}^{*,j},\\
		u_t^j = u_{t|t}^{*,j}.
	\end{aligned}
\end{equation}
Moreover, the terminal constraint (\ref{eq:Constraints5}) enforces $x^{*,j}_{t+N|t} \in \mathcal{CS}^{j-1}$ and, from (\ref{p^*}) and (\ref{eq:LMPCobj}),
\begin{equation} \label{eq:x_opt_n}
	x^{*,j}_{t+N|t} = \sum_{k=0}^{j-1} \sum_{t=0}^\infty \lambda^{*,k}_t x_t^k.
\end{equation}
We define 
\begin{equation}
	\bar{u} = \sum_{k=0}^{j-1} \sum_{t=0}^\infty \lambda^{*,k}_t u_t^k, ~ \in \mathcal{U},
\end{equation}
and
\begin{equation} \label{eq:x_defined}
	\begin{aligned}
		\bar{x} &= A x^{*,j}_{t+N|t} + B \bar{u} = \sum_{k=0}^{j-1} \sum_{t=0}^\infty \lambda_t^{*,k} \Big(A x_t^k + B u_t^k \Big) = \\
		&= \sum_{k=0}^{j-1} \sum_{t=0}^\infty \lambda_t^{*,k} x_{t+1}^k \in \mathcal{CS}^{j-1}.
	\end{aligned}
\end{equation} 
Since the state update in (\ref{eq:system}) and (\ref{eq:Constraints1}) are assumed identical we have that
\begin{equation}\label{eq:NoMissmatch}
	x_{t+1}^j = x_{t+1|t}^{*,j}.
\end{equation}
At time $t+1$ of the $j$-th iteration the input sequence and the related feasible state trajectory
\begin{subequations} \label{eq:Feasible}
	\begin{align}
		&[u_{t+1|t}^{*,j},~u_{t+2|t}^{*,j},~...,~u_{t+N-1|t}^{*,j},~\bar{u}], \label{eq:Feasible_1} \\
		&[x_{t+1|t}^{*,j},~x_{t+2|t}^{*,j},~...,~x_{t+N-1|t}^{*,j},~x_{t+N|t}^{*,j},~\bar{x}] \label{eq:Feasible_2}
	\end{align}
\end{subequations}
satisfy input and state constrains (\ref{eq:Constraints1})-(\ref{eq:Constraints2})-(\ref{eq:Constraints4}). Therefore, (\ref{eq:Feasible}) is a feasible solution for the LMPC (\ref{eq:Constraints}) and (\ref{eq:MPC}) at time $t+1$.\\
We showed that at the $j$-th iteration, $\forall j \geq 1$ , \emph{(i):} the LMPC is feasible at time $t=0$ and  \emph{(ii):} if the LMPC is feasible at time $t$, then the LMPC is feasible at time $t+1$. Thus, we conclude by induction that the LMPC in (\ref{eq:Constraints}) and (\ref{eq:MPC}) is feasible $\forall j \geq 1$ and $t \in\zz_{0+}$.

Next we use the fact the Problem (\ref{eq:Constraints}) is time-invariant at each iteration $j$ and we replace $J_{t\rightarrow t+N}^{\scalebox{0.4}{LMPC},j}(\cdot)$ with $J_{0\rightarrow N}^{\scalebox{0.4}{LMPC},j}(\cdot)$. In order to show the asymptotic stability of $x_F$ we have to show that the optimal cost, $J_{0\rightarrow N}^{\scalebox{0.4}{LMPC},j}(\cdot)$, is a Lyapunov function for the equilibrium point $x_F$ (\ref{eq:Def_x_f}) of the closed loop system (\ref{eq:system}) and (\ref{eq:MPC}) \cite{c16}. Continuity of $J_{0\rightarrow N}^{\scalebox{0.4}{LMPC},j}(\cdot)$ can be shown as in \cite{c12}. Moreover from (\ref{eq:Constraints1InfCost}), $J_{0\rightarrow N}^{\scalebox{0.4}{LMPC},j}(x) \succ 0 ~ \forall ~ x \in \rr^n \setminus \{x_F\}$ and $J_{0\rightarrow N}^{\scalebox{0.4}{LMPC},j}(x_F)=0$. Thus, we need to show that  $J_{0\rightarrow N}^{\scalebox{0.4}{LMPC},j}(\cdot)$ is decreasing along the closed loop trajectory.\\
From (\ref{eq:NoMissmatch}) we have $x_{t+1|t}^{*,j}=x_{t+1}^j$, which implies that
\begin{equation}\label{eq:LyapProof_01}
	\begin{aligned}
		J_{0\rightarrow N}^{\scalebox{0.4}{LMPC},j}(x_{t+1|t}^*) = J_{0\rightarrow N}^{\scalebox{0.4}{LMPC},j}(x_{t+1}^j).
	\end{aligned}
\end{equation}
Given the optimal input sequence and the related optimal trajectory in (\ref{eq:OptimalSolutionMPC}) and the definition of the $P^{j-1}(\cdot)$ (\ref{eq:Psol}), the optimal cost is given by
\begin{equation}\label{eq:LyapProof}
	\begin{aligned}
		J_{0\rightarrow N}^{\scalebox{0.4}{LMPC},j}(x_t^j&)=\min_{u_{t|t},\ldots,u_{t+N-1|t}} \bigg[  \sum_{k=0}^{N-1}  h(x_{k|t},u_{k|t}) +\\
		&~~~~~~~~~~~~~~~~~~~~~~~~~~~~~~~~~+ P^{j-1}(x_{N|t}) \bigg] = \\
		= h(x_{t|t}^{*,j}&,u_{t|t}^{*,j}) + \sum_{k=1}^{N-1}  h(x^{*,j}_{t+k|t},u^{*,j}_{t+k|t}) + P^{j-1}(x^{*,j}_{t+N|t}) =\\
		= h(x_{t|t}^{*,j}&,u_{t|t}^{*,j}) + \sum_{k=1}^{N-1}  h(x^{*,j}_{t+k|t},u^{*,j}_{t+k|t}) + \\
		&~~~~~~~~~~~~~~~~~~~~+\sum_{k =0}^{j-1} \sum_{t=0}^{\infty} \lambda^{*,k}_t \sum_{l=0}^{\infty} h(x_{t+l}^k, u_{t+l}^k).
	\end{aligned}
\end{equation}
We can further simplify the above expression using (\ref{eq:p_term}), (\ref{eq:x_opt_n})-(\ref{eq:x_defined}) and the fact that $h(\cdot, \cdot)$ is jointly convex in the arguments,
\begin{equation}\label{eq:LyapProof1}
	\begin{aligned}
		J_{0\rightarrow N}^{\scalebox{0.4}{LMPC},j}(x_t^j&)
		= h(x_{t|t}^{*,j},u_{t|t}^{*,j}) + \sum_{k=1}^{N-1}  h(x^{*,j}_{t+k|t},u^{*,j}_{t+k|t})+\\
		+\sum_{k =0}^{j-1} \sum_{t=0}^{\infty} & \lambda^{*,k}_t h(x_t^k, u_t^k) + \sum_{k =0}^{j-1} \sum_{t=0}^{\infty} \lambda^{*,k}_t \sum_{l=1}^{\infty} h(x_{t+l}^k, u_{t+l}^k) \\
		&\geq h(x_{t|t}^{*,j},u_{t|t}^{*,j}) + \sum_{k=1}^{N-1}  h(x^{*,j}_{t+k|t},u^{*,j}_{t+k|t}) + \\
		&+ h\Big(\sum_{k =0}^{j-1} \sum_{t=0}^{\infty} \lambda^{*,k}_t x_t^k, \sum_{k =0}^{j-1} \sum_{t=0}^{\infty} \lambda^{*,k}_t u_t^k\Big) + \\ &~~~~~~~~~~~~~~~~+ \sum_{k =0}^{j-1} \sum_{t=0}^{\infty} \lambda^{*,k}_t J_{t\rightarrow \infty}^{k}(x_{t+1}^{k})\geq \\
		\geq h(x_{t|t}^{*,j}&,u_{t|t}^{*,j}) +  \sum_{k=1}^{N-1}  h(x^{*,j}_{t+k|t},u^{*,j}_{t+k|t}) + h(x^{*,j}_{t+N|t},\bar{u}) ~ + \\
		&~~~~~~~~~~~~~~~~~~~~~~~~~~~~~~~~~~~~~~+ P^{j-1}(\bar{x}) \geq \\
		&\geq h(x_{t|t}^{*,j},u_{t|t}^{*,j}) + J_{0\rightarrow N}^{\scalebox{0.4}{LMPC},j}(x_{t+1|t}^{*,j}).
	\end{aligned}
\end{equation}
Note that, in the above derivation, we used the fact that $\bar \lambda_0^{k} = 0$ and $\bar \lambda_{t+1}^{k} =\lambda_t^{*,k},~\forall k \in \{ 0, j-1 \}, ~ t \in \zz_{0+}$ is a feasible solution to problem (\ref{p^*}) and therefore $\sum_{t=0}^{\infty} \lambda^{*,k}_t J_{t\rightarrow \infty}^{k}(x_{t+1}^{k})$ is a upper bound for $P^{j-1}(\bar{x})$. Finally, from equations (\ref{eq:MPC}), (\ref{eq:BC}) and (\ref{eq:LyapProof_01})-(\ref{eq:LyapProof1}) we conclude that the optimal cost is a decreasing Lyapunov function along the closed loop trajectory,
\begin{equation}\label{eq:LyapProof2}
	\begin{aligned}
		J_{0\rightarrow N}^{\scalebox{0.4}{LMPC},j}(x_{t+1}^j)-J_{0\rightarrow N}^{\scalebox{0.4}{LMPC},j}(x_{t}^j) \leq - &h(x_{t}^{j},u_{t}^{j}) < 0, \\
		&\forall~x_t^j \in R^n \setminus \{x_F\}
	\end{aligned}
\end{equation}
Equation (\ref{eq:LyapProof2}), the positive definitiveness of $h(\cdot, \cdot)$ and the continuity of $J_{0\rightarrow N}^{\scalebox{0.4}{LMPC},j}(\cdot)$ imply that $x_F$ is asymptotically stable. \cvd

\subsection{Convergence properties}
In this Section we assume that the LMPC (\ref{eq:Constraints}) and (\ref{eq:MPC}) converges to a steady state trajectory. We show two results. First, the $j$-th iteration cost $J_{0\rightarrow \infty}^{j}(\cdot)$ does not worsen as $j$ increases. Second, the steady state trajectory is the solution to the infinite horizon control problem (\ref{eq:ConstraintsInf}).

\begin{theorem}
	\label{th2}
	Consider system (\ref{eq:system}) in closed loop with the LMPC controller (\ref{eq:Constraints}) and (\ref{eq:MPC}).
	Let $\mathcal{CS}^j$ be the convex safe set at the $j$-th iteration as defined in (\ref{eq:CS}). Let assumption 1 hold, then the iteration cost $J_{0\rightarrow \infty}^{j}(\cdot)$ does not increase with the iteration index $j$.
\end{theorem}

\textbf{Proof:}
Follows from \textit{Theorem 2} in \cite{LMPC} \cvd 

 \begin{theorem}
 	\label{th3}
 	Consider system (\ref{eq:system}) in closed loop with the LMPC controller (\ref{eq:Constraints}) and (\ref{eq:MPC}) with $N > 1$.
 	Let $\mathcal{CS}^j$ be the convex safe set at the $j$-th iteration as defined in (\ref{eq:CS}). Let assumption 1 hold and assume that the closed loop system (\ref{eq:system}) and (\ref{eq:MPC}) converges to a steady state trajectory ${\bf{x}}^\infty$, for iteration $j\rightarrow \infty$. Then, the steady state input ${\bf{u}}^\infty = \lim_{j \to \infty} {\bf{u}}^j$ and the related steady state trajectory ${\bf{x}}^\infty = \lim_{j \to \infty} {\bf{x}}^j$ is a global optimal solution for the infinite horizon optimal control problem (\ref{eq:ConstraintsInf}), i.e., ${\bf{x}}^\infty = {\bf{x}}^{*}$ and ${\bf{u}}^\infty = {\bf{u}}^{*}$.
 \end{theorem}
 \begin{theorem}
 	Consider system (\ref{eq:system}) in closed loop with the LMPC controller (\ref{eq:Constraints}) and (\ref{eq:MPC}) with $N > 1$.
 	Let $\mathcal{CS}^j$ be the sampled safe set at the $j$th iteration as defined in (\ref{eq:CS}). Let assumption 1 hold and assume that the closed loop system (\ref{eq:system}) and (\ref{eq:MPC}) converges to a steady state trajectory ${\bf{x}}^\infty$, for iteration $j\rightarrow \infty$. Denote $\Int ( \mathcal{S})$ as the interior of the set $\mathcal{S}$, and recall the definition of ones-step predecessor $\Pre(\cdot)$ and successor  $\Succ(\cdot)$ sets from \cite[Section II]{LMPC}.
 	If $x_k^{\infty} \in \Int ( \Pre (x_{k+1}^{\infty}))$ and $x_{k+1}^{\infty} \in \Int( \Succ(x_{k}^{\infty}))$ for all $k \geq 0$. Then, the steady state input ${\bf{u}}^\infty = \lim_{j \to \infty} {\bf{u}}^j$ and the related steady state trajectory ${\bf{x}}^\infty = \lim_{j \to \infty} {\bf{x}}^j$ is a global optimal solution for the infinite horizon optimal control problem (\ref{eq:ConstraintsInf}), i.e., ${\bf{x}}^\infty = {\bf{x}}^{*}$ and ${\bf{u}}^\infty = {\bf{u}}^{*}$.
 \end{theorem}

\textbf{Proof:} Follows from the convexity of (\ref{eq:Constraints}), (\ref{eq:MPC}) and of Problem (\ref{eq:ConstraintsInf}) and  \textit{Theorem 3} in \cite{LMPC}.

\cvd

\section{Example: Constrained LQR controller}
In this section, we test the proposed LMPC for linear system on the following infinite horizon linear quadratic regulator with constraints (CLQR)
\begin{subequations}\label{eq:CLQR}
	\begin{align}
		J_{0\rightarrow \infty}^*(x_S)&=\min_{u_0, u_1,\ldots} \sum\limits_{k=0}^{\infty} \Big[ ||x_k||_2^2 + ||u_k||_2^2 \Big] \\
		\textrm{s.t. }
		&x_{k+1}= \begin{bmatrix} 1 & 1 \\ 0 & 1 \end{bmatrix} x_k +  \begin{bmatrix} 0 \\ 1 \end{bmatrix} u_k,~\forall k\geq 0 \label{eq:CLQR1}\\
		&x_0=x_S,\label{eq:CLQR2}\\
		& \begin{bmatrix} -4 \\ -4 \end{bmatrix} \leq x_k \leq \begin{bmatrix} 4 \\ 4 \end{bmatrix} ~ \forall k\geq 0\label{eq:CLQR3} \\
		&-1 \leq u_k \leq 1 ~~\forall k\geq 0. \label{eq:CLQR4}
	\end{align}
\end{subequations}

In \cite{LMPC} we showed that the LMPC converges to the solution of the infinite horizon control problem (\ref{eq:CLQR}), whenever we have can compute a feasible trajectory ${\bf{x^0}}$. However, the LMPC in \cite{LMPC} is implemented using the sampled Safe Set (\ref{eq:SS}) as a terminal constraint, instead of the proposed convex Safe Set (\ref{eq:CS}). Therefore, also for linear systems, the LMPC presented in \cite{LMPC} involves the solution of a Mixed Integer Programming (MIP) program which is computationally expensive. In the following, we show that the proposed convex formulation of the LMPC for linear systems reduces the computational burden by several order of magnitude and it converges to the solution of the infinite horizon control problem (\ref{eq:CLQR}). 

The LMPC (\ref{eq:Constraints}), (\ref{eq:MPC}) is implemented with the quadratic running cost $h(x_k,u_k) = ||x_k||_2^2 + ||u_k||_2^2$, an horizon of $\bar{N}$ steps, and the states and input constraints (\ref{eq:CLQR3})-(\ref{eq:CLQR4}). The LMPC (\ref{eq:Constraints}) and (\ref{eq:MPC}) is reformulated as a Quadratic Programming and it is implemented in YALMIP (\cite{c28}) using the solver quadprog. In order to implement the terminal cost (\ref{eq:p^*Rel}) we defined the time ${t}_j$ at which the iterations is completed,
\begin{equation}\label{eq:Terminate}
	\begin{aligned}
		{t}_j = \min \textrm{} \Big\{ t\in\zz_{0+}  : J_{0\rightarrow \infty}^{\scalebox{0.4}{LMPC},j}(x_t^j) \leq \epsilon \Big\}.
	\end{aligned}
\end{equation}
with $\epsilon = 10^{-8}$.

For $x_S = [-3.95, -0.05]^T$ and $\bar{N}=4$, the LMPC converges a to steady state solution ${\bf{x}^{\infty}}={\bf{x}^{8}}$, after $8$ iterations, with a error of $\gamma= 10^{-10}$:
\begin{equation}\label{eq:Steady}
	\begin{aligned}
		\max_{t \in [0, {t}_9]} ||{{x}}^9_t-{{x}}^8_t||_2 < \gamma
	\end{aligned}.
\end{equation}

Table \ref{table:Cost} shows the evolution of the iterations cost. We notice that accordingly with \textit{Theorem 2} the cost is non-increasing over the iterations.
\textbf{\begin{table}[h!]
\caption{Optimal cost of the LMCPC at each $j$-th iteration}\label{table:Cost}
\centering
\begin{tabular}{lll}
\multicolumn{3}{l}{~~\textbf{Iteration}}{~~\textbf{Iteration Cost}}                   \\ \hline
 & $j = 0$    & ~~~~$57.1959612235$ \\
 & $j =1$     & ~~~~$49.9313760793$  \\
 & $j =2$     & ~~~~$49.9166091658$  \\
 & $j =3$     & ~~~~$49.9163668042$  \\
 & $j =4$     & ~~~~$49.9163602472$  \\
 & $j =5$     & ~~~~$49.9163600537$  \\
 & $j =6$     & ~~~~$49.9163600469$  \\
 & $j =7$     & ~~~~$49.9163600468$  \\
 & $j =8$     & ~~~~$49.9163600464$ 
\end{tabular}
\end{table}}
	
	Furthermore, the solution of the LMPC for linear system is compared with the exact solution of the CLQR (\ref{eq:CLQR}), which is computed using the algorithm in \cite{c16}. Given the optimal solution to the infinite horizon optimal control problem (\ref{eq:CLQR}),
	\begin{equation}\label{eq:Optimaltask}
		\begin{aligned}
			{\bf{x}}^{*} ~ = ~ [{{x}}^{*}_0,~{{x}}^{*}_1,~...,~{{x}}^{*}_t,~...],&&\\
			{\bf{u}}^{*} ~ = ~ [{{u}}^{*}_0,~{{u}}^{*}_1,~...,~{{u}}^{*}_t,~...],&&\\
		\end{aligned}
	\end{equation}
	we define the approximation error as
	\begin{equation}
		\begin{gathered}\label{eq:deviation}
			\sigma_t = ||x_t^{\infty}-x^*_t||_2.
		\end{gathered}
	\end{equation}
	$\sigma_t$ quantifies, at each time step $t$, the distance between the optimal trajectory of the CLQR (\ref{eq:CLQR}) and steady state trajectory at of the LMPC (\ref{eq:Constraints}) and (\ref{eq:MPC}). The maximum approximation error is
	\begin{equation}
		\bar{\sigma} = \max[\sigma_0, \dots, \sigma_{{t}_{\infty}}] = 8.6 \times 10^{-6}.
	\end{equation}
	Moreover, the $2$-norm of the normalized difference between the exact optimal cost and the cost associated with the steady state trajectory is 
	\begin{equation}
		\Delta J =\frac{||J_{0\rightarrow \infty}^*({{x}}_S)-J_{0\rightarrow \infty}^{*,\infty}(x_0^{\infty})||_2}{J_{0\rightarrow \infty}^*({{x}}_S)} \times 100 = 1.8 \times 10^{-20}.
	\end{equation}
	The LMPC for linear system (\ref{eq:Constraints}) and (\ref{eq:MPC}) has converged to global optimal solution.

	\balance 
	
 We tested the LMPC (\ref{eq:Constraints}), (\ref{eq:MPC}) with different initial conditions $x_S$ and horizon length $\bar{N} > 1$ to experimentally validate \textit{Theorems 1-3}. Table \ref{table:Convergence} shows the maximum approximation error, $\bar{\sigma}$, and $\Delta J$. We underline that for all the tested scenarios, regardless of the horizon length, the proposed LMPC converged to the global optimal solution of the infinite horizon control problem. It is interesting to notice that the LMPC (\ref{eq:Constraints}), (\ref{eq:MPC}) with a longer horizon $\bar{N}$ has more freedom to explore the state space and therefore it converges faster to the steady state trajectory.
	\begin{table}[h!]
		\centering
		\caption{Convergence of the LMPC for different initial conditions}
		\label{table:Convergence}
		\begin{tabular}{llllc}
			{\bf{$x_S$}}            & $\bar{N}$ & $\bar{\sigma}$                   & $\Delta J$                & Iterations  \\ \hline
			$[-3.95, -0.05]^T$ &       2        & $2.6 \times 10^{-1}$         & $1.3 \times 10^{-1}$          & $44$  \\
			$[-3.95, -0.05]^T$ &       3        & $1.9 \times 10^{-5}$         & $1.7 \times 10^{-17}$         & $26$ \\
			$[-3.95, -0.05]^T$ &       4        & $8.6 \times 10^{-6}$         & $1.8 \times 10^{-20}$         & $~8$    \\
			$[-4, 0]^T$            &       2        & $3.6 \times 10^{-1}$         & $4.2 \times 10^{-1}$      & $74$    \\
			$[-4, 0]^T$            &       3        & $1.6 \times 10^{-5}$         & $5.9 \times 10^{-18}$     & $26$   \\
			$[-4, 0]^T$            &       4        & $5.2 \times 10^{-6}$         & $1.2 \times 10^{-20}$     & $~8$   \\
			$[   -2,         2]^T$ &       3        & $7.8 \times 10^{-2}$         & $5.1 \times 10^{-3}$      & $80$   \\
			$[   -2,         2]^T$ &       3        & $1.7 \times 10^{-5}$         & $3.4 \times 10^{-17}$      & $22$   \\
			$[   -2,         2]^T$ &       4        & $7.3 \times 10^{-6}$         & $5.7 \times 10^{-20}$     & $~8$   \\
			$[    0,       1.5]^T$&       2        & $1.0 \times 10^{-1}$         & $1.5 \times 10^{-2}$       & $45$   \\
			$[    0,       1.5]^T$&       3        & $1.8 \times 10^{-5}$         & $3.4 \times 10^{-17}$       & $27$   \\
			$[    0,       1.5]^T$&       4        & $6.6 \times 10^{-6}$         & $2.0 \times 10^{-19}$      & $~8$  \\
			\end{tabular}
	\end{table}

	Finally, we compare the computational burden associate with the LMPC (\ref{eq:Constraints}), (\ref{eq:MPC}) and with the LMPC in \cite{LMPC}. The proposed LMPC (\ref{eq:Constraints}), (\ref{eq:MPC}) applied to Problem (\ref{eq:CLQR}) converged in $40 s$ to a steady state trajectory. On the other hand, the LMPC in \cite{LMPC} applied to Problem (\ref{eq:CLQR})  took $2 hr$ to reach convergence. Therefore, we conclude that the proposed approach significantly reduces the computational burden of the control logic preserving the properties of the LMPC.  
	
	\section{Conclusions}
	In this paper, an extension to the learning Model Predictive Control (LMPC) is presented. The controller is designed for linear system and it significantly reduces the computational burden associated with the LMPC. A convex safe set and a terminal cost, learnt from previous iterations, allow to guarantee the recursive feasibility and stability of the closed loop system. Furthermore, the LMPC is guaranteed to improve the performance of the close-loop system over the iterations. We tested the proposed control logic on an infinite horizon linear quadratic regulator with constraints (CLQR) to show that the proposed control logic converges to the optimal solution of the infinite optimal control problem. Finally, we compared the computation time of the proposed strategy with the computational time of the LMPC for nonlinear system, and we showed that the proposed control logic reduces the computational burden by several order of magnitudes.

\bibliography{ifacconf}             

\begin{thebibliography}{15}
\providecommand{\natexlab}[1]{#1}
\providecommand{\url}[1]{\texttt{#1}}
\providecommand{\urlprefix}{URL }
\expandafter\ifx\csname urlstyle\endcsname\relax
  \providecommand{\doi}[1]{doi:\discretionary{}{}{}#1}\else
  \providecommand{\doi}{doi:\discretionary{}{}{}\begingroup
  \urlstyle{rm}\Url}\fi

\bibitem[{Borrelli(2003)}]{c16}
Borrelli, F. (2003).
\newblock \emph{Constrained optimal control of linear and hybrid systems},
  volume 290.
\newblock Springer.

\bibitem[{Bristow et~al.(2006)Bristow, Tharayil, and Alleyne}]{c3}
Bristow, D.A., Tharayil, M., and Alleyne, A.G. (2006).
\newblock A survey of iterative learning control.
\newblock \emph{IEEE Control Systems}, 26(3), 96--114.

\bibitem[{Garcia et~al.(1989)Garcia, Prett, and Morari}]{c11}
Garcia, C.E., Prett, D.M., and Morari, M. (1989).
\newblock Model predictive control: theory and practice-a survey.
\newblock \emph{Automatica}, 25(3), 335--348.

\bibitem[{Jones and Morari(2010)}]{jones2010polytopic}
Jones, C.N. and Morari, M. (2010).
\newblock Polytopic approximation of explicit model predictive controllers.
\newblock \emph{IEEE Transactions on Automatic Control}, 55(11), 2542--2553.

\bibitem[{Lee and Lee(2007)}]{c7}
Lee, J.H. and Lee, K.S. (2007).
\newblock Iterative learning control applied to batch processes: An overview.
\newblock \emph{Control Engineering Practice}, 15(10), 1306--1318.

\bibitem[{Lee et~al.(2000)Lee, Lee, and Kim}]{lee2000model}
Lee, J.H., Lee, K.S., and Kim, W.C. (2000).
\newblock Model-based iterative learning control with a quadratic criterion for
  time-varying linear systems.
\newblock \emph{Automatica}, 36(5), 641--657.

\bibitem[{Lee and Lee(2000)}]{lee2000convergence}
Lee, K.S. and Lee, J.H. (2000).
\newblock Convergence of constrained model-based predictive control for batch
  processes.
\newblock \emph{IEEE Transactions on Automatic Control}, 45(10), 1928--1932.

\bibitem[{Lofberg(2004)}]{c28}
Lofberg, J. (2004).
\newblock Yalmip: A toolbox for modeling and optimization in matlab.
\newblock In \emph{Computer Aided Control Systems Design, 2004 IEEE
  International Symposium on}, 284--289. IEEE.

\bibitem[{Mayne et~al.(2000)Mayne, Rawlings, Rao, and Scokaert}]{c12}
Mayne, D.Q., Rawlings, J.B., Rao, C.V., and Scokaert, P.O. (2000).
\newblock Constrained model predictive control: Stability and optimality.
\newblock \emph{Automatica}, 36(6), 789--814.

\bibitem[{Rosolia and Borrelli(2017)}]{LMPC}
Rosolia, U. and Borrelli, F. (2017).
\newblock Learning model predictive control for iterative tasks. a data-driven
  control framework.
\newblock \emph{IEEE Transactions on Automatic Control}, 63(7), 1883--1896.

\bibitem[{Rucco et~al.(2015)Rucco, Notarstefano, and Hauser}]{c18}
Rucco, A., Notarstefano, G., and Hauser, J. (2015).
\newblock An efficient minimum-time trajectory generation strategy for
  two-track car vehicles.
\newblock \emph{IEEE Transactions on Control Systems Technology}, 23(4),
  1505--1519.

\bibitem[{Sharp and Peng(2011)}]{c17}
Sharp, R. and Peng, H. (2011).
\newblock Vehicle dynamics applications of optimal control theory.
\newblock \emph{Vehicle System Dynamics}, 49(7), 1073--1111.

\bibitem[{Subbaraman and Benosman(2016)}]{subbaraman2016extremum}
Subbaraman, A. and Benosman, M. (2016).
\newblock Extremum seeking-based iterative learning model predictive control
  (esilc-mpc).
\newblock \emph{IFAC-PapersOnLine}, 49(13), 193--198.

\bibitem[{Tamar et~al.(2016)Tamar, Thomas, Zhang, Levine, and
  Abbeel}]{2016learning}
Tamar, A., Thomas, G., Zhang, T., Levine, S., and Abbeel, P. (2016).
\newblock Learning from the hindsight plan--episodic mpc improvement.
\newblock \emph{arXiv preprint arXiv:1609.09001}.

\bibitem[{Wang et~al.(2009)Wang, Gao, and Doyle}]{c8}
Wang, Y., Gao, F., and Doyle, F.J. (2009).
\newblock Survey on iterative learning control, repetitive control, and
  run-to-run control.
\newblock \emph{Journal of Process Control}, 19(10), 1589--1600.

\end{thebibliography}

                                                   







\appendix
\end{document}